\documentclass[conference]{IEEEtran}
\usepackage{stmaryrd}
\usepackage{amsfonts}
\usepackage{bm,algorithmicx,algpseudocode,threeparttable,slashbox,booktabs}

\usepackage{graphicx,times,amsmath} 

\IEEEoverridecommandlockouts    

\begin{document}

\title{\ \\ \LARGE\bf A Comparative Study of STA on Large Scale Global Optimization \thanks{$^{*}$Corresponding author of this paper. The authors are with the School of Information Science and Engineering, Central South University, Changsha 410083, China (email: ychh@csu.edu.cn).} \thanks{This work was supported by the National Natural Science
Foundation of China (Grant No. 61503416, 61533020, 61533021,61590921).}}

\author{Xiaojun Zhou, Chunhua Yang$^{*}$ and Weihua Gui}


\maketitle

\begin{abstract}
State transition algorithm has been emerging as a new intelligent global optimization method in recent few years.
The standard continuous STA has demonstrated powerful
global search ability for global optimization problems whose dimension is no more than 100.
In this study, we give a test report to present the performance of standard continuous STA
for large scale global optimization when compared with other state-of-the-art evolutionary algorithms. From the experimental results, it is shown that the standard continuous STA
still works well for almost all of the test problems, and its global search ability
is much superior to its competitors.
\end{abstract}


\section{Introduction}

\PARstart{S}{tate transition algorithm} (STA) has been emerging as a new
intelligent optimization method for global optimization in recent few years \cite{Zhou2011initial}-\cite{Zhou2015optimal}. In state transition algorithm, a solution to
an optimization problem is considered as a state, and an update of a solution can be
regarded as a state transition.
By referring to state space representation, on the basis of current state $\bm x_k$, the unified form of generation of a new state $\bm x_{k+1}$ in state transition algorithm can be described as follows:
\begin{equation}
\left \{ \begin{array}{ll}
\bm x_{k+1}= A_{k} \bm x_{k} + B_{k} \bm u_{k}\\
y_{k+1}= f(\bm x_{k+1})
\end{array} \right.,
\end{equation}
where $\bm x_{k} = [x_1, x_2, \cdots, x_n]^T$ stands for a state, corresponding to a solution of an optimization problem; $\bm u_{k}$ is a function of $\bm x_{k}$ and historical states;
$A_{k}$ and
$B_{k}$ are state transition matrices, which are usually some state transformation operators;
$f(\cdot)$ is the objective function or fitness function, and
$y_{k+1}$ is the function value at $\bm x_{k+1}$.

Unlike most of the existing evolutionary algorithms, the basic STA
is an individual-based iterative method. Based on an current state, a regular neighborhood is automatically
formed by using certain state transformation operators, since there exists stochastic properties
in the state transition matrices, and then a sampling technique is used to create a candidate state set.
That is to say, the generation of a candidate set in STA is completely different from most other evolutionary algorithms.
Furthermore, special local and global search operators are both designed, and in the meanwhile,
there exists an alternative way of using local and global operators in STA.
The form of STA can be continuous or discrete, called continuous STA or discrete STA respectively,  depending on the
state transformation operators. In continuous STA, we have designed four state transformation operators named \emph{rotation}, \emph{translation}, \emph{expansion}, and \emph{axesion} to deal with continuous variables (see \cite{Zhou2012} for details); while in discrete STA, other four state transformation operators named \emph{swap}, \emph{shift}, \emph{symmetry} and
\emph{substitute} are designed as well, and they can tackle discrete variable in an effective way (please refer to \cite{Zhou2015discrete} for details). The powerfulness of both continuous and discrete STA has been demonstrated in
\cite{Zhou2012}-\cite{Wang2016sta} in terms of global search ability and convergence rate. In this study, we focus on the continuous STA for the following global optimization problem
\begin{eqnarray}
\min_{\bm x \in \Omega} f(\bm x)
\end{eqnarray}
where $\bm x \in \mathbb{R}^n$, $\Omega \subseteq \mathbb{R}^n$ is a closed and compact set, which is usually composed of lower and upper bounds of $\bm x$.

The effectiveness and efficiency of continuous STA have been testified when compared with other state-of-the-art
intelligent optimization methods, like real-coded genetic algorithm (RCGA) \cite{Tran2010real}, comprehensive learning particle swarm optimizer (CLPSO) \cite{Liang2006comprehensive}, self-adaptive differential evolution (SaDE) \cite{Qin2009differential} and artificial bee colony (ABC) algorithm \cite{Karaboga}.
However, in these studies, the size of the benchmark functions chosen for test is no more than 100.
It is reported that the performance of most intelligent optimization methods will deteriorate severely when
the size increases, especially for large scale global optimization.
Therefore, the motivation of this study is to test the continuous STA on large scale global optimization problems.

The remainder of this paper is organized as follows:
Section II gives a brier review of continuous STA and its procedures.
Section III presents the experimental results and discussions of continuous STA with its competitors on large scale benchmark problems.
The conclusions and future perspectives are given in Section IV.
\section{A brief review of continuous STA}
The initial version of continuous STA was firstly proposed in \cite{Zhou2011initial}, in which, there are only three state transformation operators, and then in \cite{Zhou2011new}, the axesion transformation was replenished to strengthen single dimensional search. By replenishing the axesion transformation and changing the rotation factor $\alpha$ to  decline periodically in an outer loop, the standard continuous STA was born in \cite{Zhou2012}.
\subsection{State transition operators}
Using state space transformation for reference, four special
state transformation operators are designed to generate continuous solutions for an optimization problem.\\
(1) Rotation transformation
\begin{equation}
\bm x_{k+1}= \bm x_{k}+\alpha \frac{1}{n \|\bm x_{k}\|_{2}} R_{r} \bm x_{k},
\end{equation}
where $\alpha$ is a positive constant, called the rotation factor;
$R_{r}$ $\in$ $\mathbb{R}^{n\times n}$, is a random matrix with its entries being uniformly distributed random variables defined on the interval [-1, 1],
and $\|\cdot\|_{2}$ is the 2-norm of a vector. This rotation transformation
has the function of searching in a hypersphere with the maximal radius $\alpha$. \\
(2) Translation transformation\\
\begin{equation}
\bm x_{k+1} = \bm x_{k}+  \beta  R_{t}  \frac{\bm x_{k}- \bm x_{k-1}}{\|\bm x_{k}- \bm x_{k-1}\|_{2}},
\end{equation}
where $\beta$ is a positive constant, called the translation factor; $R_{t}$ $\in \mathbb{R}$ is a uniformly distributed random variable defined on the interval [0,1].
The translation transformation has the function of searching along a line from $x_{k-1}$ to $x_{k}$ at the starting point $x_{k}$ with the maximum length $\beta$.
\\
(3) Expansion transformation\\
\begin{equation}
\bm x_{k+1} = \bm x_{k}+  \gamma  R_{e} \bm x_{k},
\end{equation}
where $\gamma$ is a positive constant, called the expansion factor; $R_{e} \in \mathbb{R}^{n \times n}$ is a random diagonal
matrix with its entries obeying the Gaussian distribution. The expansion transformation
has the function of expanding the entries in $\bm x_{k}$ to the range of [-$\infty$, +$\infty$], searching in the whole space.\\
(4) Axesion transformation\\
\begin{equation}
\bm x_{k+1} = \bm x_{k}+  \delta  R_{a}  \bm x_{k}\\
\end{equation}
where $\delta$ is a positive constant, called the axesion factor; $R_{a}$ $\in \mathbb{R}^{n \times n}$ is a random diagonal matrix with its entries obeying the Gaussian distribution and only one random position having nonzero value. The axesion transformation is aiming to search along the axes, strengthening single dimensional search.
\subsection{Regular neighborhood and sampling}
For a given solution, a candidate solution is generated by using one of the aforementioned state transition operators.
Since the state transition matrix in each state transformation is random, the generated candidate solution is not unique.
Based on the same given point, it is not difficult to imagine that a regular neighborhood will be automatically formed
when using certain state transition operator.
In theory, the number of candidate solutions in the neighborhood is infinity;
as a result, it is impractical to enumerate all possible candidate solutions.

Since the entries in state transition matrix obey certain stochastic
distribution, for any given solution, the new candidate becomes a random vector and its corresponding solution
(the value of a random vector) can be regarded as a sample.
Considering that any two random state transition matrices in each state transformation
are independent, several times of state transformation (called the degree of search enforcement, \textit{SE} for short)
based on the same given solution are performed
for certain state transition operator, consisting of \textit{SE} samples.
It is not difficult to find that all of the \textit{SE} samples are independent, and they are
representatives of the neighborhood.
Taking the rotation transformation for example, a total number of
\textit{SE} samples are generated in pseudocode as follows
\begin{algorithmic}[1]
\For{$i\gets 1$, \textit{SE}}
\State $\mathrm{State}(:,i) \gets {\mathrm{Best}} +\alpha \frac{1}{n \| {\mathrm{Best}} \|_{2}} R_{r}  {\mathrm{Best}} $
\EndFor
\end{algorithmic}
where $Best$ is the incumbent best solution, and \textit{SE} samples are stored in
the matrix $\mathrm{State}$.

\subsection{An update strategy}
As mentioned above, based on the incumbent best solution, a total number of
\textit{SE} candidate solutions are generated.
A new best solution is selected from the candidate set by virtue of the fitness function, denoted as \emph{newBest}.
Then, an update strategy based on greedy criterion is used to update the incumbent best as shown below
\begin{equation}
Best =
\left\{ \begin{aligned}
&\mathrm{newBest},  \;\;\;\mathrm{if}\; f(\mathrm{newBest}) < f(\mathrm{Best}) \\
&\mathrm{Best},\;\;\;\;\;\;    \;\;\;\mathrm{otherwise}
\end{aligned} \right.
\end{equation}
\subsection{Algorithm procedure of the basic continuous STA}
With the state transformation operators, sampling technique and update strategy, the basic state transition algorithm can be described by the following pseudocode
\begin{algorithmic}[1]
\Repeat
    \If{$\alpha < \alpha_{\min}$}
    \State {$\alpha \gets \alpha_{\max}$}
    \EndIf
    \State {Best $\gets$ expansion(funfcn,Best,SE,$\beta$,$\gamma$)}
    \State {Best $\gets$ rotation(funfcn,Best,SE,$\alpha$,$\beta$)}
    \State {Best $\gets$ axesion(funfcn,Best,SE,$\beta$,$\delta$)}
    \State {$\alpha \gets \frac{\alpha}{\textit{fc}}$}
\Until{the specified termination criterion is met}
\end{algorithmic}

\indent As for detailed explanations, rotation$(\cdot)$ in above pseudocode is given for illustration purposes as follows
\begin{algorithmic}[1]
\State{oldBest $\gets$ Best}
\State{fBest $\gets$ feval(funfcn,oldBest)}
\State{State $\gets$ op\_rotate(Best,SE,$\alpha$)}
\State{[newBest,fnewBest] $\gets$ fitness(funfcn,State)}
\If{fnewBest $<$ fBest}
    \State{fBest $\gets$ fnewBest}
    \State{Best $\gets$ newBest}
    \State{State $\gets$ op\_translate(oldBest,newBest,SE,$\beta$)}
    \State{[newBest,fnewBest] $\gets$ fitness(funfcn,State)}
    \If{fnewBest $<$ fBest}
        \State{fBest $\gets$ fnewBest}
        \State{Best $\gets$ newBest}
    \EndIf
\EndIf
\end{algorithmic}

As shown in the above pseudocodes, the rotation factor $\alpha$ is decreasing periodically from a maximum
value $\alpha_{\max}$ to a minimum value $\alpha_{\min}$ in an
exponential way with base \textit{fc}, which is called lessening coefficient.
op\_rotate$(\cdot)$ and op\_translate$(\cdot)$ represent the implementations of proposed sampling technique for rotation and
translation operators, respectively, and fitness$(\cdot)$ represents the implementation of selecting the new best solution
from \textit{SE} samples. It should be noted that the translation operator is only executed
when a solution better than the incumbent best solution can be found in the \textit{SE} samples
from rotation, expansion or axesion transformation.
In the basic continuous STA, the parameter settings are given as follows:
$\alpha_{\max} = 1, \alpha_{\min} = 1e$-4, $\beta = 1, \gamma = 1, \delta = 1$, $\textit{SE} = 30, \textit{fc} = 2$.

When using the fitness$(\cdot)$ function, solutions in \emph{State} are projected into
$\Omega$ by using the following formula
\begin{equation}
x_i =
\left\{ \begin{aligned}
&u_i,  \;\;\; \mathrm{if}\; x_i > u_i \\
&l_i,  \;\;\; \mathrm{if}\; x_i < l_i\\
&x_i,  \;\;\; \mathrm{otherwise}
\end{aligned} \right.
\end{equation}
where $u_i$ and $l_i$ are the upper and lower bounds of $x_i$ respectively.

\section{Experimental test for large scale global optimization problems}
In the experiment, the standard continuous STA is coded in MATLAB R2010b
on Intel(R) Core(TM) i5-4200U CPU @1.60GHz under
Window 7 environment. The maximum number of iterations (MaxIter for short) is chosen to
terminate the program. More specifically, the MaxIter for 100D, 200D and 500D problem is 1e3, 2e3, 5e3 respectively, and
a total of 10 independent runs are performed with randomly chosen initial points.
In the same time, the comprehensive learning particle swarm optimizer (CLPSO) \cite{Liang2006comprehensive}, self-adaptive differential evolution (SaDE) \cite{Qin2009differential} are used for comparison
with the same parameter settings.

\subsection{Benchmark functions}
Five well-known benchmark functions are used for test as follows

(1) Spherical function
\[
f_1= \sum_{i=1}^n x_{i}^{2},
\]
where the global optimum $\bm x^{*} = (0, \cdots, 0)$ and $f(\bm x^{*}) = 0$,
$-100 \leq x_i \leq 100, i= 1, \cdots, n$.

(2) Rosenbrock function
\[
f_2= \sum_{i=1}^n (100(x_{i+1}-x_{i}^2)^2 + (x_{i}-1)^2),
\]
where the global optimum $\bm x^{*} = (1, \cdots, 1)$ and $f(\bm x^{*}) = 0$,
$-30 \leq x_i \leq 30, i= 1, \cdots, n$.

(3) Rastrigin function
\[
f_3 =  \sum_{i=1}^n(x_{i}^{2}-10\cos(2 \pi x_{i})+10),
\]
where the global optimum $\bm x^{*} = (0, \cdots, 0)$ and $f(\bm x^{*}) = 0$,
$-5.12 \leq x_i \leq 5.12, i= 1, \cdots, n$.

(4) Griewank function
\[
f_4 = \frac{1}{4000} \sum_{i=1}^n x_{i}^{2}-  \prod_{i}^{n} \cos|\frac{x_{i}}{\sqrt{i}}| + 1,
\]
where the global optimum $\bm x^{*} = (0, \cdots, 0)$ and $f(\bm x^{*}) = 0$,
$-600 \leq x_i \leq 600, i= 1, \cdots, n$.

(5) Ackley function
\[
f_5(\bm x)\!=\!20\!+\!e\!-\!20\exp(\!-\!0.2\sqrt{\frac{1}{n} \sum_{i=1}^n x_{i}^2})
\!-\!\exp(\frac{1}{n} \sum_{i=1}^n \cos(2\pi x_{i}))
\]
where the global optimum $\bm x^{*} = (0, \cdots, 0)$ and $f(\bm x^{*}) = 0$,
$-32 \leq x_i \leq 32, i= 1, \cdots, n$.

\subsection{Results and Discussion}
\begin{figure}[htp]
\centerline{\includegraphics[width=3.8in]{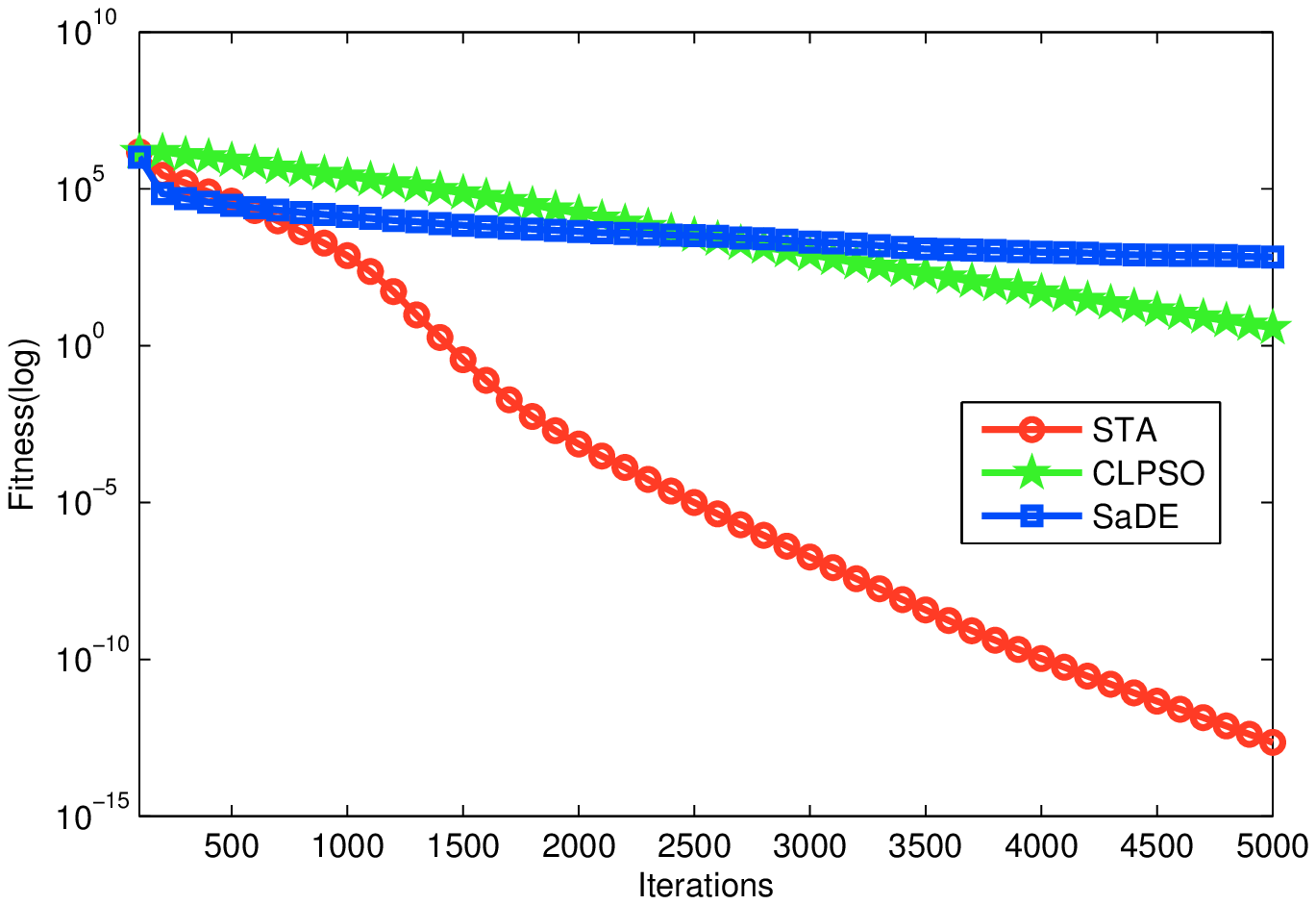}}
\caption{The average iterative curves of three algorithms on $f_1$(500D)} \label{fig_f1}
\end{figure}

\begin{figure}[htp]
\centerline{\includegraphics[width=3.8in]{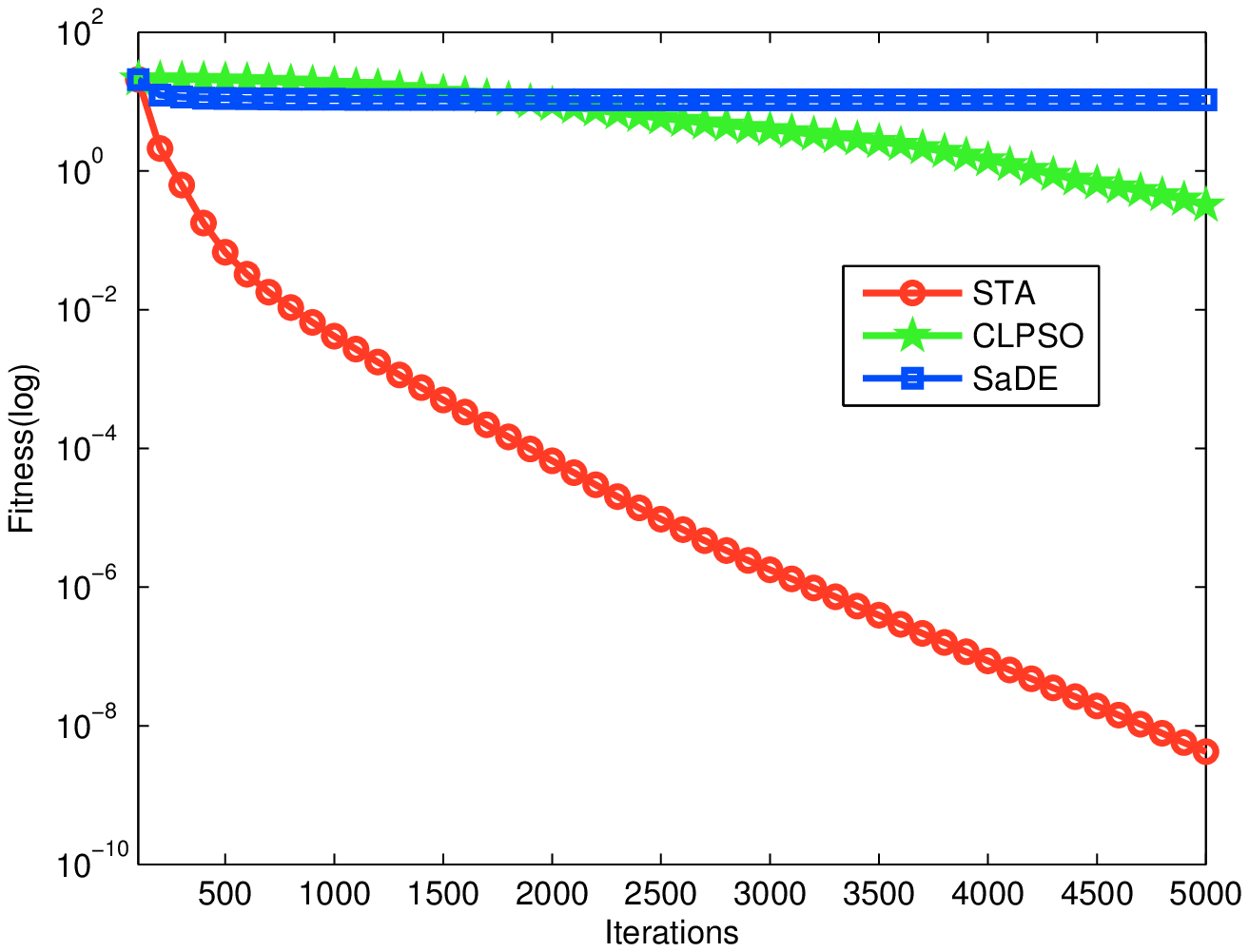}}
\caption{The average iterative curves of three algorithms on $f_5$(500D)} \label{fig_f5}
\end{figure}

Experimental results are shown in Table 1.
First and foremost, it can be found that only the standard continuous STA
can find global solutions or approximate global solutions with very high precision for
almost all test functions from 100 dimension to 500 dimension except the Rosenbrock function.
In the same time,
SaDE can find approximate global solutions for the Sphere function only for 100 dimension, and fails for other benchmark functions within such a prescribed iteration number, and the CLPSO fails for all the benchmark function test in such a situation.
Then, for the Rastrigin function and Griewank function with 100 dimension, the standard
continuous STA can find their global solutions with extremely high precision. While for Spherical function and Ackley function, the standard continuous STA has the potential to find their global solutions with very high precision as well, as indicated in
the iterative curves of Fig.\ref{fig_f1} and Fig.\ref{fig_f5}.

\begin{figure}[htp]
\centerline{\includegraphics[width=3.8in]{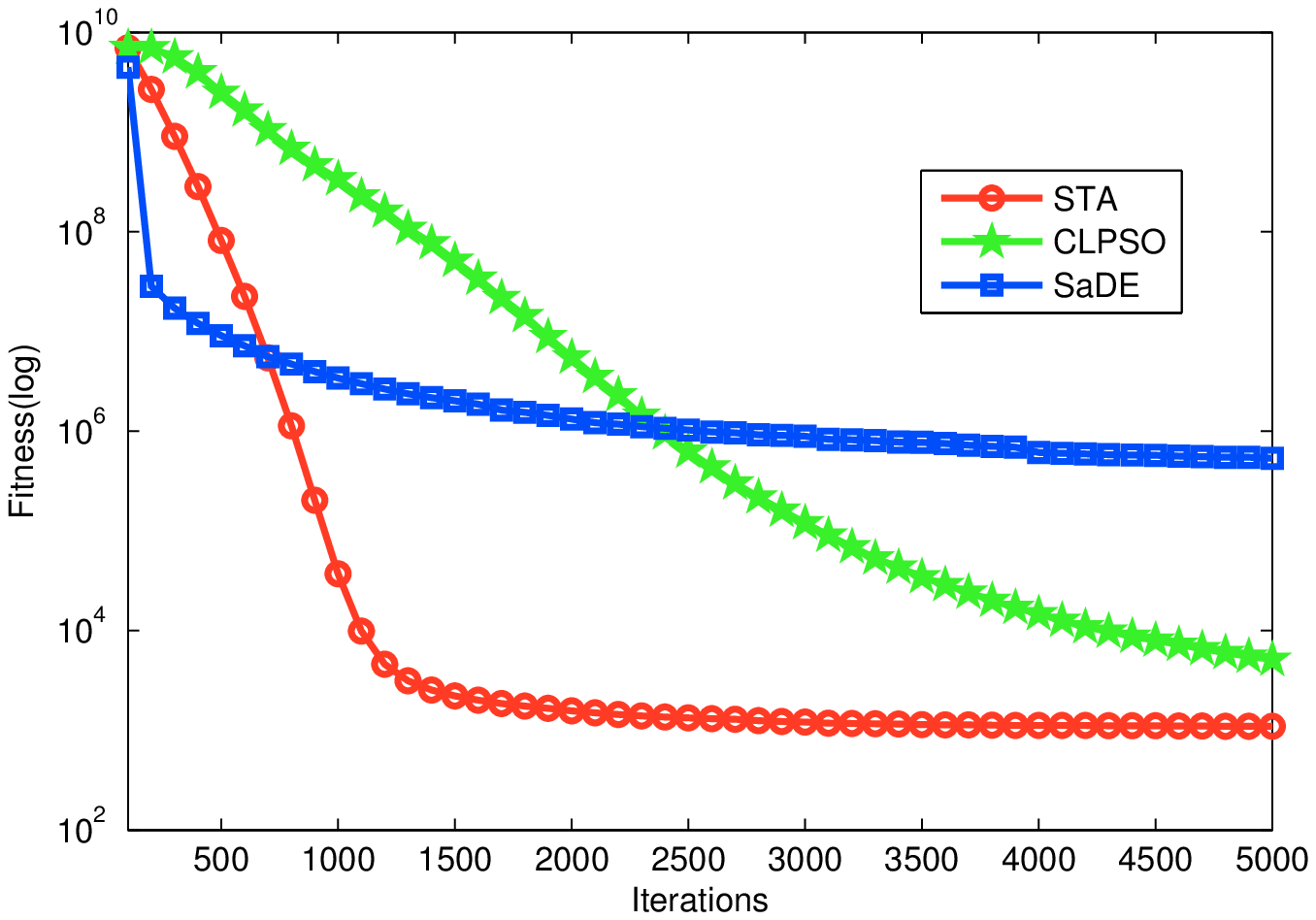}}
\caption{The average iterative curves of three algorithms on $f_2$(500D)} \label{fig_f2}
\end{figure}

\begin{table*}[!htbp]
\caption{Comparative results of continuous STA for large scale benchmark functions}
\centering
\footnotesize
\begin{tabular}{{ccccc}}
\hline
\toprule[1pt]
Function & Algorithm &  100D & 200D & 500D   \\
\hline
        & STA	  & \textbf{8.6933e-122 $\pm$ 2.7150e-121}  &  \textbf{4.0937e-16 $\pm$  1.4755e-16} & \textbf{1.2176e-13 $\pm$  1.9606e-13}\\
$f_{1}$ & CLPSO  & 4.5374 $\pm$ 0.7757 &  4.0776 $\pm$ 1.3061  & 2.8750 $\pm$ 0.4522\\
        & SaDE    & 1.4385e-7 $\pm$ 1.6453e-7 &  0.0757 $\pm$ 0.1033  & 623.1392 $\pm$ 411.6156 \\
\hline
        & STA	  & 1\textbf{66.4412 $\pm$ 54.4169}  &	\textbf{398.7591  $\pm$   127.213}8       & \textbf{1.0925e3  $\pm$    116.3676}        \\
$f_{2}$ & CLPSO  & 3.3276e3  $\pm$ 504.8178 &	3.4533e3  $\pm$   387.6409       & 4.8055e3  $\pm$  274.4664          \\
        & SaDE    & 351.3534 $\pm$ 49.1156  &	1.1426e3  $\pm$   326.3298	     & 2.1565e5	 $\pm$  1.8998e5              \\
\hline
        & STA	  & \textbf{0  $\pm$ 0}             &	\textbf{4.4611e-11 $\pm$ 7.9241e-11}          & \textbf{9.1495e-11  $\pm$  4.6790e-11}        \\
$f_{3}$ & CLPSO  & 115.0501 $\pm$ 13.4231  &	255.5583  $\pm$ 20.1691         & 670.4771 $\pm$ 23.1019            \\
        & SaDE    & 34.5702 $\pm$ 9.4130   &	111.8367 $\pm$ 14.5737	        & 398.0244 $\pm$ 31.5910               \\
\hline
        & STA	  & \textbf{0  $\pm$  0 }      &	\textbf{2.1094e-16 $\pm$ 3.5108e-17}         & \textbf{3.2629e-14  $\pm$ 9.9673e-14}         \\
$f_{4}$ & CLPSO  & 0.9945 $\pm$ 0.0594           &	0.7268 $\pm$ 0.0598          & 0.2783  $\pm$ 0.0177              \\
       & SaDE    & 0.0414 $\pm$ 0.0536           &	0.1670 $\pm$ 0.2450	        & 4.8419 $\pm$ 4.1143               \\
\hline
        & STA	  & \textbf{3.0198e-15 $\pm$ 1.1235e-15}         &	\textbf{3.5181e-9 $\pm$  1.9309e-9}        & \textbf{3.1285e-9  $\pm$ 5.1074e-10}        \\
$f_{5}$ & CLPSO  & 1.4243  $\pm$ 0.2634          &	0.6200 $\pm$ 0.1205          & 0.2729 $\pm$ 0.0560              \\
        & SaDE    &  4.2657  $\pm$ 0.7385         &	7.2220 $\pm$ 0.6434	        & 10.5653 $\pm$ 0.6280              \\
\bottomrule[1pt]
\hline
\end{tabular}
\end{table*}

\begin{figure}[htp]
\centerline{\includegraphics[width=3.8in]{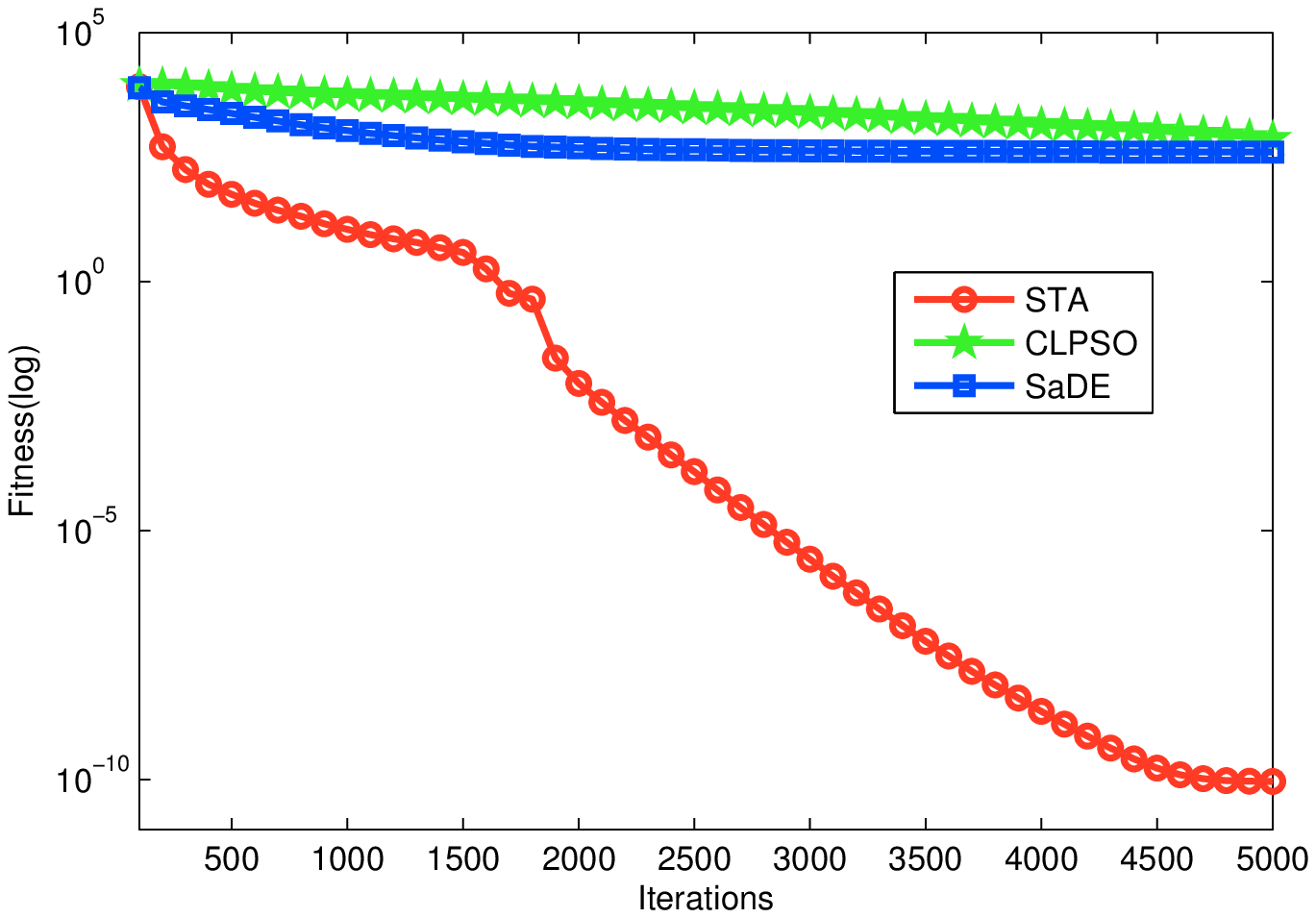}}
\caption{The average iterative curves of three algorithms on $f_3$(500D)} \label{fig_f3}
\end{figure}

\begin{figure}[htp]
\centerline{\includegraphics[width=3.8in]{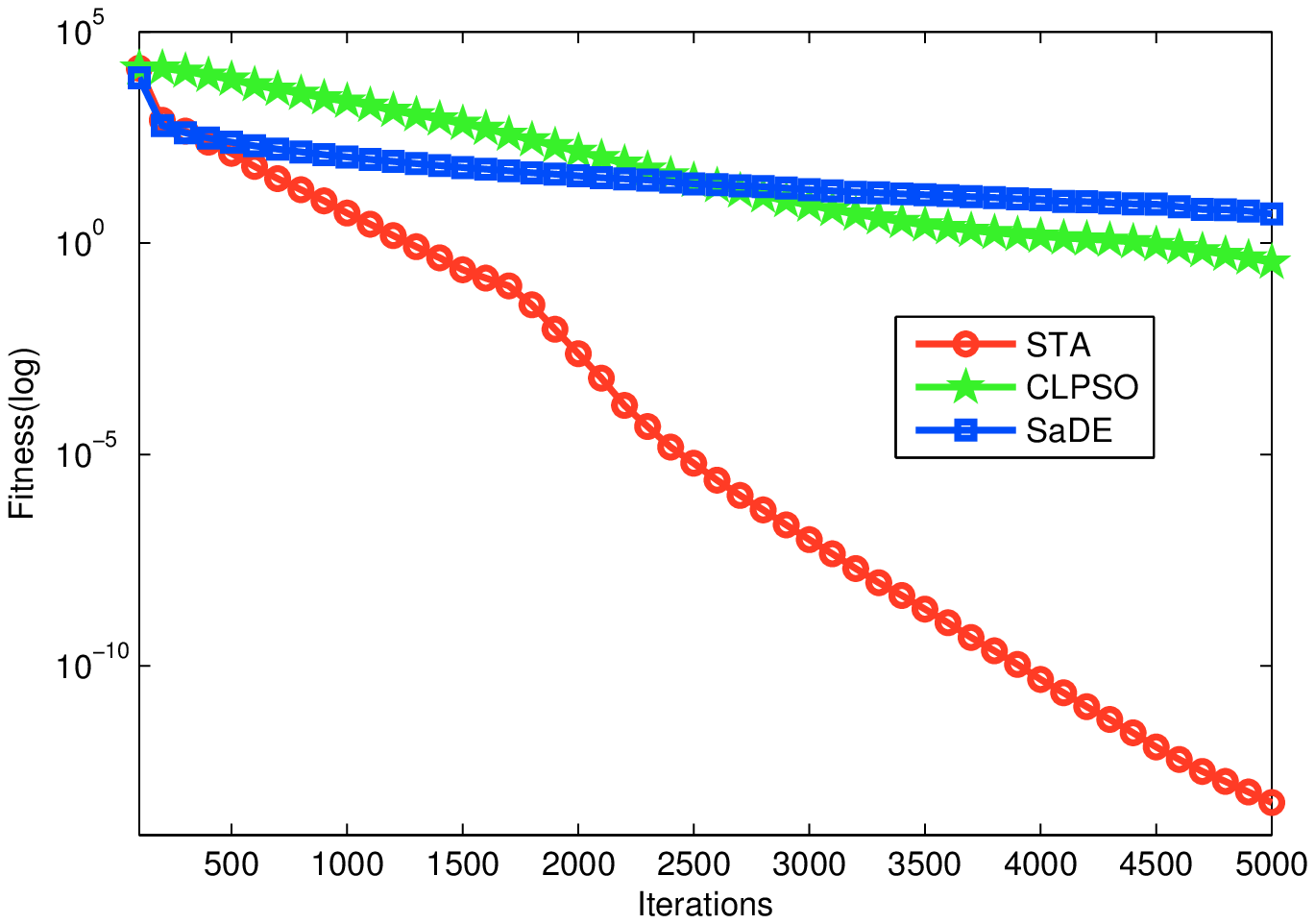}}
\caption{The average iterative curves of three algorithms on $f_4$(500D)} \label{fig_f4}
\end{figure}

It was reported that CLPSO and SaDE performed very nice for optimization problems with low dimension (see \cite{Liang2006comprehensive} and \cite{Qin2009differential}) .
However, within the prescribed few iterations, their performance deteriorates sharply for large scale optimization problems (500D), as shown from Fig.\ref{fig_f1} to Fig.\ref{fig_f5}, which indicates that they have very slow convergence rate for large scale optimization problems.
On the other hand, it should be noted that the standard continuous STA performs not well for the Rosenbrock function because the iterative curves decrease very slowly at the late stage.
These observed phenomena have shown that the standard continuous STA has very strong global search ability while
the local search ability should be strengthened, since for the Rosenbrock function, the local search is more important in the later search stage.

\section{Conclusion and future perspectives}

A comparative study of the standard continuous STA for large scale global optimization is reported in this paper.
Experimental results have demonstrated that the standard continuous STA has very strong global search ability.
However, it should be noted that local search ability of standard continuous STA is not so good.
In our future work, we will focus on the local search of continuous STA to accelerate its convergence rate.

On the other hand, the presented experimental results are based on integration.
In recent few decades, the divide-and-conquer strategy has been widely applied to large-scale global optimization via decomposition (see \cite{Yang2008} and \cite{Li2012}).
In our future work, this strategy will be adopted into STA to decompose the original optimization problem to
a series of sub-optimization problems and the corresponding composition techniques will be studied as well.



%

\def\V{\rm vol.~}
\def\N{no.~}
\def\pp{pp.~}
\def\Pot{\it Proc. }
\def\IJCNN{\it International Joint Conference on Neural Networks\rm }
\def\ACC{\it American Control Conference\rm }
\def\SMC{\it IEEE Trans. Systems\rm , \it Man\rm , and \it Cybernetics\rm }

\def\handb{ \it Handbook of Intelligent Control: Neural\rm , \it
    Fuzzy\rm , \it and Adaptive Approaches \rm }

\end{document}